\documentclass[a4paper,11pt]{article}
\usepackage{amscd,amssymb,amsthm}
\usepackage{amssymb,palatino}
\usepackage{euscript}
\usepackage{color}

\newtheorem{theorem}{Theorem}

\begin{document}

\title{Continued Fractions, Quadratic Fields, and Factoring: \\
        Some Computational Aspects         
 } %
\author{Michele Elia
\thanks{Politecnico di Torino
Corso Duca degli Abruzzi 24, I - 10129 Torino -- Italy; ~~ e-mail: elia@polito.it }}


\maketitle

\thispagestyle{empty}

\begin{abstract}
\noindent
 Legendre discovered  that the continued fraction expansion of $\sqrt N$ having odd period leads directly to an explicit representation of $N$ as the sum of two squares. 
In this vein, it was recently observed that the continued fraction expansion of $\sqrt N$ having even period directly produces a factor of composite $N$.
It is proved here that these apparently fortuitous occurrences allow us to propose and apply a variation
of Shanks' infrastructural method which significantly reduces the asymptotic computational burden with respect to currently used factoring techniques.
\end{abstract}

\section{Introduction}
In a letter to { Pierre de Carcavi}, August $14^{th}$ $1659$, { Pierre de Fermat} 
reported several propositions; in particular, he stated the following theorem: 
{\it  Every prime $p$ of the form $4k+1$ is uniquely expressible as the sum of two squares, i.e.}
$  p = X^2+Y^2 ~~  \Leftrightarrow ~~  p \equiv 1  \bmod 4$,  ~
whose first known proof was given by Euler
using Fermat's {\em infinite descent} method. Many other proofs have been given, some
 constructive, others non-constructive; in particular, among the latter, 
Zagier's one-sentence proof deserves to be mentioned for its conciseness \cite{zagier}.
Among the numerous constructive proofs, 
two different proofs by Gauss stand out.
The first is direct, and gives  
$   x= \frac{(2k)!}{2(k!)^2}   \bmod p$  and  $y= \frac{((2k)!)^2}{2(k!)^2}   \bmod p $;
 the partially incomplete proof was completed, a century later, by Jacobsthal.  
The second proof is based on quadratic forms of discriminant $-4$, and considers two equivalent
principal quadratic forms with discriminant $-4$:
$    p X^2+2 b_1XY+ \frac{b_1^2+1}{p}Y^2 ~~\mbox{and}~~~  x^2+y^2$,
where $b_1$ is a root of $z^2+1$ modulo $p$.
The first form represents $p$ trivially with $X=1$ and $Y=0$, thus
Gauss' reduction produces the unique reduced form
 in the class \cite{mathews}, and meanwhile yields $x$ and $y$.

\noindent
 Jacobsthal's constructive solution  (1906) is based on counting the number of points on the elliptic curve
$y^2=n(n^2-a)$ in $\mathbb Z_p$.  He considers the sum of Legendre symbols
$$  S(a) = \sum_{n=1}^{p-1}  \left(\frac{n(n^2-a)}{p} \right)  \Rightarrow x=\frac{1}{2} S(q_R) ~~,~~ 
   y=\frac{1}{2} S(q_N) $$
 where $q_R, q_N \in \mathbb Z_p$ are any quadratic residue and non-residue, respectively,
 \cite{jacobsthal}.

\noindent
Legendre's proof is reported on pages 59-60 of \cite{legendre}. 
It is constructive, since it yields  $X$ and $Y$
 from the complete remainder of the continued fraction expansion of $\sqrt p$.
It is well explained in his own words

\begin{quotation}
{\em ... Donc tous le fois que l'\'equation  $x^2 -Ay^2=-1$ est r\'esoluble (ce qui ha lieu entre autre cas
 lorsque $A$ est un numbre premier $4n+1$) le nombre $A$ peut toujours \^etre decompos\'e en deux
 quarr\'es; et cette d\'ecomposition est donn\'ee immediatement par lo quotient-complet
 $\frac{\sqrt A + I}{D}$  qui r\'epond au second des quotients moyens compris dans la premi\`ere p\'eriode du d\'eveloppement de
$\sqrt A$; le nombres $I$ et $D$  \'etant ainsi connu, on aura $A=D^2+I^2$. 

\vspace{2mm}
Cette conclusion ranferme un des plus beaux th\'eor\`emes de la science des nombres, savoir,
{\em que tout nombre premier $4n+1$ est la somme de deux quarr\'es;}  elle donne en m\^eme temps 
  le moyen de faire cette d\'ecomposition d'une mani\`ere directe et sans aucun}
{\bf \textcolor{black}{t\^atonnement}}.
\end{quotation}

\noindent
Thus, Legendre's proof gives the representation of any composite $N$ such that the period of the continued fraction for
 $\sqrt N$ is odd, or equivalently, $x^2-N y^2=-1$ is solvable in integers \cite{legendre,sierp,elia1}. \\
As a counterpart to Legendre's finding,  when the period of the continued fraction
 expansion of $\sqrt N$ is even, we directly obtain, under mild conditions, a factor of a composite $N$.
In particular, this is certainly the case when both prime factors of $N=pq$ are congruent $3$ modulo $4$ \cite{elia1}. 
Legendre's solution of Fermat's theorem tacitly introduces a connection between
 continued fractions and  the ramified primes of quadratic number fields, obviously without using
 this notion more than a century before Dedekind's invention. 
To explain these singular connections, the paper is organized as follows. 
Section 2 summarizes the properties of the continued fraction expansion of $\sqrt N$.  
Section 3 discusses the factorization of composite numbers $N$.
Lastly, Section 4  draws conclusions.

\section{Preliminaries}
A regular continued fraction is an expression of the form
\begin{equation}
  \label{cf}
    a_0+\frac{1}{ a_1+\frac{1}{a_2+\frac{1}{ a_3+ \cdots}}} ~~,
\end{equation}
 where $a_0$, $a_1$,  $a_2, \ldots, a_i, \ldots$ is a sequence, possibly infinite,
 of positive integers. 
A convergent of a continued fraction is the sequence of fractions $\frac{A_m}{B_m}$,
 each of which is obtained by truncating the continued fraction at the  $(m+1)$-th term. 
The fraction  $\frac{A_m}{B_m}$ is called the $m$-th convergent \cite{dave,hardy}. 
A continued fraction is said to be definitively periodic, 
 with period $\tau$, if, starting from a finite position $n_o$, a fixed pattern
 $a_1'$,  $a_2', \ldots, a_{\tau}'$ repeats indefinitely. 
Lagrange showed that any definitively periodic continued fraction, of period length $\tau$, represents
 a positive number of the form $a+b\sqrt{N}$, $a,b \in \mathbb Q$, i.e. an element of 
$\mathbb Q(\sqrt N)$, and conversely
 any such positive number is represented by a definitively periodic continued fraction
 \cite{dave,sierp}. 
The period of the continued fraction expansion of $\sqrt{N}$ begins
 immediately after the first term $a_0$, and is written as
$         \sqrt{N} = \left[ a_0 , \overline{a_1,a_2, \ldots, a_2, a_1, 2a_0} \right]$,
where the over-lined part is the period, which includes a
 palindromic part formed by the $\tau-1$ terms $~~a_1,a_2, \ldots, a_2, a_1$. 
In Carr's book \cite[p.70-71]{carr} we find a good collection
 of properties of the continued fraction expansion of $\sqrt{N}$, which are summarized in the
 following, along with some properties taken from \cite{dave,sierp}.

\begin{enumerate}
 \item Let $c_n$ and $r_n$ be the elements of two sequences of positive integers defined by
   the relation
$$  \frac{\sqrt{N}+c_n}{r_n}=a_{n+1}+\frac{r_{n+1}}{\sqrt{N}+c_{n+1}} $$
with $c_0= \left\lfloor  \sqrt N \right\rfloor$, and $ r_0 =N-a_0^2$;
 the elements of the sequence $a_1, a_2, \ldots , a_n \ldots$
   are thus obtained as the integer parts of the left-side fraction,
  which is known as the complete quotient.
 \item Let $a_0= \lfloor \sqrt{N} \rfloor$ be initially computed, and set $c_0 = a_0$,  $r_0 =N-a_0^2$,
 then sequences $\{ c_n \}_{n \geq 0}$  and $\{ r_n \}_{n \geq 0}$ are produced by the 
recursions
{\small
\begin{equation}
   \label{contfrac}
    a_{m+1} = \left\lfloor \frac{a_0+c_m}{r_m} \right\rfloor ~,~ c_{m+1}=a_{m+1} r_m -c_m 
     ~,~ r_{m+1} = \frac{N-c_{m+1}^2}{r_m}.
\end{equation}
}
These recursions allow us to compute the sequence
 $\{a_m\}_{m\geq 1}$ using only rational arithmetical operations, and
 the iterations may be stopped when $a_m = 2 a_0$, having completed a period.
\item  If the period length $\tau$ is odd, set $\ell= \frac{\tau-1}{2}$; Legendre discovered
 and proved that the complete quotient $\frac{\sqrt{N} +c_\ell}{r_\ell}$ gives a representation of 
 $N=c_\ell^2+r_\ell^2$ as the sum of two squares.
 \item Numerator $A_n$ and denominator $B_n$ of the $n$-th convergent to $\sqrt{N}$
  can be recursively computed as  $A_n = a_n A_{n-1}+ A_{n-2}$
   and $B_n=a_n B_{n-1}+ B_{n-2}$, ~$n \geq $1,
 respectively, with initial conditions $A_{-1}=1$, $B_{-1}=0$, $A_{0}=a_0$, and $B_{0}=1$.
The numerator $A_m$ and the denominator $B_m$ of any convergent are shown to be
 relatively prime by the relation  
 $ A_mB_{m-1}-A_{m-1}B_m=(-1)^{m-1} $ \cite[p.85]{dave}.
\item Using the sequences $\{ A_m \}_{m \geq 0}$ and $\{ B_m \}_{m \geq 0}$, two sequences 
$\mathbf \Delta = \{\Delta_m=A_m^2-N B_m^2\}_{m \geq 0}$, and
$\mathbf \Omega = \{ \Omega_m = A_m A_{m-1}-N B_m B_{m-1} \}_{m \geq 1}$ are introduced.
It can easily be checked that $\Omega_m^2 - \Delta_m \Delta_{m-1} =N,~ \forall m \geq 1$. 
 The elements of $\mathbf \Delta$ and  $\mathbf \Omega$ satisfy a system of linear recurrences 
   \begin{equation}
         \label{Deltarecur}
         \left\{
           \begin{array}{ll}
             \Delta_{m+1} =  a_{m+1}^2 \Delta_{m}  + 2 a_{m+1} \Omega_{m}+  \Delta_{m-1} \\
             \Omega_{m+1} =  \Omega_{m} + a_{m+1} \Delta_{m}  \\
           \end{array}  \right.  ~~~~~~m \geq 1 
      \end{equation}
    with initial conditions $\Delta_0=a_0^2-N$, $\Delta_1=(1+a_0a_1)^2-N a_1^2$ and
    $\Omega_1=(1+a_0a_1) a_0-N a_1$. By (\ref{Deltarecur}), it is immediate to see that
    $c_{m+1}= |\Omega_m|$ and $r_{m+1}=|\Delta_m|$.
\item The period of  $\mathbf \Delta$ and  $\mathbf \Omega$ is $\tau$ or $2\tau$, depending on
  whether $\tau$ is even or odd.
\item The sequence of ratios $ \frac{A_n}{B_n}$ assumes the limit value $\sqrt{N}$ as $n$ 
 goes to infinity, due to the inequality
$  \left|   \frac{A_n}{B_n} -\sqrt N \right|  \leq   \frac{1}{B_n B_{n+1}}  ~~,  $
  since $A_n$ and $B_n$ go to infinity along with $n$. Since $\frac{A_n}{B_n} <\sqrt N$, 
  if $n$ is even, and  $\frac{A_n}{B_n} >\sqrt N$, if $n$ is odd \cite{hardy},
 any convergent of even index is smaller than any convergent of odd index. This property implies that 
 the terms of the sequence $\mathbf \Delta$ have 
 alternating signs, with $\Delta_1 > 0$.
\item The value $c_0=a_0$ is the greatest value that $c_n$ may assume.  
           No $a_n$ or $r_n$ can be greater than $2a_0$. \\
           If $r_n=1$ then $a_{n+1}=a_0$. 
           For all $n$ greater than $0$, we have $a_0-c_n < r_n\leq 2a_0$. 
           The first complete quotient that is repeated 
          is $\frac{\sqrt{N} +c_0}{r_0}$, and $a_1$, $r_0$, and $c_0$ commence
          each cycle of repeated terms.
 \item Through the first period, we have the equalities
    $a_{\tau-j}=a_j$ , $r_{\tau-j-2}=r_j$, and $c_{\tau-j-1}=c_j$. 
 \item The period $\tau$ has the tight upper bound
$      0.72 \sqrt{N} \ln N ~,~ \  N > 7$, 
     as was shown by Kraitchik \cite[p.95]{steuding}.
   However, the period length has irregular behavior as a function of $N$, because it may assume any
   value from $1$, when $N=M^2+1$, 
     to values close to the order $O( \sqrt{N} \ln N )$  \cite{sierp}.
  \item  Define the sequence of quadratic forms 
           $\mathbf f_m(x,y)=\Delta_m x^2+ 2 \Omega_m x y + \Delta_{m-1}y^2$, $m \geq 1$, 
           which has the same period as $\mathbf \Delta$.
           Every $\mathbf f_m(x,y)$ is a reduced form of discriminant $4N$. 
          Within the first block, all quadratic forms $ \mathbf f_m(x,y)$,  $1 \leq m \leq \tau$  are distinct,
          and constitute the principal class $\mathbf \Gamma(\mathbf f)$  of reduced forms,
         with the  ordering of the elements inherited from $\mathbf \Delta$.
      The definition of reduced form used here is slightly different from the classic one: set  
  $\kappa=\min \{ |\Delta_m |, |\Delta_{m-1} | \}$; it is easily checked that $\Omega_m$ is the sole
 integer such that $\sqrt N- |\Omega_m|<\kappa <\sqrt N+|\Omega_m| $, with the sign of $\Omega_m$
 chosen opposite to the sign of $\Delta_m $. 
         Since the sign of $\Delta_{m-1}$ is the same as that of $\Omega_m$, which  is opposite to that of $\Delta_m$,
          in $\mathbf{\Gamma}(f)$ the two triples of signs (signatures) $(-,+,+)$ and $(+,-,-)$ alternate. 
\end{enumerate}

\noindent
The following theorems are taken, without proof,  from \cite{elia1}.

\begin{theorem}   
    \label{per2}
Starting with $m=1$, the sequences $\mathbf \Delta = \{ \Delta_m \}_{m \geq 0}$ and
 $\mathbf \Omega = \{ \Omega_m \}_{m \geq 0}$ are periodic with the same period $\tau$
 or $2\tau$ depending on whether $\tau$ is even or odd. 
The elements of the blocks  
$\{ \Delta_m \}_{m=0}^{\tau}$
 and 
 $\{ \Omega_m \}_{m=1}^{\tau}$
 satisfy the symmetry relations
 $\Delta_m=(-1)^{\tau}\Delta_{\tau-m-2},  ~~\forall~m \leq \tau-3$ and 
 $\Omega_{\tau-m-1}=(-1)^{\tau+1}\Omega_{m}, ~~\forall~m \leq \tau-2$, respectively.
\end{theorem}

\noindent 
If $\tau$ is odd, the ordered set $\{ \Delta_m \}_{m=1}^{\tau}$ has a central term of
   index $\ell= \frac{\tau-1}{2}$, $\Delta_\ell =- \Delta_{\ell-1}$ since $\tau-\ell-2=\ell-1$,
   and the equation $\Omega_\ell^2 - \Delta_\ell \Delta_{\ell-1} =N$ gives a solution of the Diophantine
   equation $x^2+y^2=N$ with $x=\Delta_\ell$ and  $y=\Omega_\ell$, the situation first recognized by Legendre. 

\noindent
If $\tau$ is even, the ordered set $\{ \Delta_m \}_{m=1}^{\tau}$ has no central term;
in this case, with  $\ell = \frac{\tau-2}{2}$ we have $\Omega_{\ell+1}  =- \Omega_\ell$
 and $\Delta_{\ell+1} = \Delta_\ell$, hence $\mathbf f_{\ell+1}(x,y)=\mathbf f_\ell(y,-x)$. 

\begin{theorem} 
    \label{sym3a}
Let the period $\tau$ of the continued fraction expansion of $\sqrt{N}$ be even; we have
 $\Omega_{\tau-1}=-a_0$,
$ \Delta_{\tau} = \Delta_{\tau-2}$, and $\Omega_{\tau}=- \Omega_{\tau-1}$.
Defining the integer $\gamma \in \mathfrak O_{\mathbb Q(\sqrt N)}$ by the product
$$  \gamma= \prod_{m=1}^\tau \left(\sqrt N +(-1)^{m} \Omega_m \right) ~~,$$
 let $\sigma$ denote the Galois automorphism of $\mathbb Q(\sqrt N)$ (i.e. $\sigma(\sqrt N=-\sqrt N$), then  $ \frac{\gamma}{\sigma(\gamma)}=A_{\tau-1}+B_{\tau-1}\sqrt{N}$ is a positive fundamental unit 
  (or the cube of the fundamental unit) of $\mathbb Q(\sqrt N)$.
\end{theorem}

\noindent
Based on this theorem, we say that the unit $\mathfrak c_{\tau-1}=A_{\tau-1}+B_{\tau-1}\sqrt{N}$ in $\mathbb Q(\sqrt N)$)
 splits $N$, if $N_1=\gcd \{A_{\tau-1} -1, N \}$ is neither $1$ nor $N$.
 Then we have the proper factorization $N=N_1 N_2$. Further, using the following involutory matrix,
   \cite{elia1}, whose square is $(-1)^{\tau} I_2$
$$    M_{\tau-1}    =  \left[
      \begin{array}{cc}
       - A_{\tau-1}    &  N B_{\tau-1} \\
         - B_{\tau-1}  &  A_{\tau-1}
     \end{array}      \right] ~~,
$$ 
 it is shown that
\begin{equation}
    \label{involutmat}
   A_{\tau-m-2} =(-1)^{m-1}( A_{\tau-1} A_m - N B_{\tau-1} B_m)  ~~1\leq m \leq \tau-2  ~~.  
\end{equation}
As an immediate consequence of this equation, if the unit  $\mathfrak c_{\tau-1}$ splits $N$, then
any pair$(A_m,A_{\tau-m-2})$ splits $N$, since taking $A_{\tau-m-2}$
modulo $N$ we have  $A_{\tau-m-2} = (-1)^{m-1} A_m  A_{\tau-1} \bmod N$ ,
thus $A_{\tau-m-2}$  is certainly different from $A_m$, because $A_{\tau-1} \neq \pm 1 \bmod N$.

\begin{theorem}
   \label{locfactor}
If the period $\tau$ of the continued fraction expansion of $\sqrt{N}$ is even, the element
 $\mathfrak c_{\tau-1}$ in $\mathbb Q(\sqrt N)$ splits $4N$, and a factor of $4N$ is located at 
 positions $\frac{\tau-2}{2}+j\tau$, $j=0,1, \ldots$, in the sequence
 $\mathbf \Delta=\{\mathfrak c_{m} \sigma(\mathfrak c_{m}) \}_{m \geq 1}$.
\end{theorem}

\section{Factorization}
Gauss recognized that the factoring problem was to be important, although very difficult,
\begin{quotation}
\noindent
{\em $\ldots$  Problema, numeros primos a compositis dignoscendi, hosque in factores
   suos primos resolvendi, ad gravissima ac utilissima totius arithmeticae pertinere, et 
   geometrarum tum veterum tum recentiorum industriam ac sagacitatem occupavisse, 
   tam notum est, ut de hac re copiose loqui superfluum foret.  $\ldots$ }
\hfill {\scriptsize \sc C. F. Gauss [{\em Disquisitiones Arithmeticae} Art. 329]}
\end{quotation}
In spite of much effort, various different approaches, and the increased importance stemming from the large
 number of cryptographic applications, no satisfactorily factoring method has yet been found. 
However, approaches to factoring based on continued fractions have lead to some of the most efficient
 factoring algorithms. 
In the following, a new variant of Shanks' infrastructural method \cite{shanks} is described
 which exploits the property of the block
 $\mathbf \Delta_1=\{ \Delta_m \}_{m=1}^{\tau}$,
 which is made more precise in the following theorem taken without proof from \cite{elia1}.

\begin{theorem}
    \label{mainfactor}
Let $N$ be a positive square-free integer. If the norm of the positive fundamental unit $\mathfrak u \in  \mathbb Q(\sqrt N)$
 is $1$, and some factor of $N$ is a square of a principal integral ideal in $  \mathbb Q(\sqrt N)$, then $\mathfrak u$
  is split for $N$. A proper factor of $N$ is found in position 
$\frac{\tau-2}{2}$ of $\mathbf \Delta_1$.
\end{theorem}

\noindent
It should be noted that $\mathbf \Delta_1$ offers several different ways for factoring a composite number $N$:
\begin{enumerate}
   \item If $\tau$ is even and $2$ is not a quadratic residue modulo $N$, then in position $\frac{\tau-2}{2}$
    of the sequence $\mathbf \Delta_1$ we find a factor of $N$.
   \item If $\tau$ is odd, then by Legendre's results we find a representation $N=X^2+Y^2$, which implies
    that $s_1= \frac{X}{Y} \bmod N$ is a square root of $-1$. If we are able to find another square root $s_2$  
      different from $-\frac{X}{Y} \bmod N$ (we have four different
     square roots of a quadratic residue modulo $N=p q$), then the difference $s_1-s_2$ contains a proper factor of $N$.
   \item If some square $d_o^2$ is found in the sequence $\mathbf \Delta_1$,
     it implies the equation $A_m^2-NB_m^2=d_o^2$, thus there is a chance that some proper
     factor of $N$ divides $(A_m-d_o)$ or $(A_m+d_o)$. \\
   The number of squares in $\mathbf \Delta_1$ is $O(\sqrt \tau)$, and about ($\frac{1}{2}$) of 
    these squares  factor $N$.  This method was introduced by Shanks.
    \item If equal terms $\Delta_m=\Delta_n, ~~m\neq n$  occur in $\mathbf \Delta_1$, with
       $m,n <\frac{\tau}{2}$, then  $A_m^2-A_n^2=0 \bmod N$ allows us to find two factors of $N$ by
      computing $\gcd\{A_m-A_n,N \}$ and $\gcd\{A_m+A_n,N \}$. This is an implementation of an old idea
     of Fermat's.
\end{enumerate}

\subsection{Computational issues}
By Theorem \ref{mainfactor} we know that a factor of $N$ is $\Delta_{\frac{\tau-2}{2}}$,
 which can be directly computed from the continued fraction of $\sqrt N$ in $\frac{\tau-2}{2}$ steps. 
 Unfortunately,  this number is usually prohibitively large. 
 However, if $\tau$ is known, using the baby-step/giant-step artifice, the number of steps can be reduced
to the order $O(\log_2 \tau)$. To this end, we can move through the principal class 
$\mathbf \Gamma(\mathbf f)$, of ordered quadratic forms $\mathbf f_{m}(x,y)$,
 by introducing a notion of distance between pairs of quadratic forms
 compliant with Gauss' quadratic form composition. 
The distance between two adjacent quadratic forms 
 $\mathbf f_{m+1}(x,y), \mathbf f_{m}(x,y) \in \mathbf \Gamma(\mathbf f)$
 is defined as
\begin{equation}
   \label{defdist} 
    d(\mathbf f_{m+1}, \mathbf f_{m}) =\frac{1}{2}  \ln \left( \frac{\sqrt{N}+(-1)^m\Omega_m }{\sqrt{N}-(-1)^m\Omega_m} \right)   ~~,  
\end{equation}
 and the distance between two quadratic forms $\mathbf f_{m}(x,y)$ and  $\mathbf f_n(x,y)$,
 with $m > n$, is defined as the sum
$    d(\mathbf f_{m}, \mathbf f_{n}) = \sum_{j=n}^{m-1}  d(\mathbf f_{j+1}, \mathbf f_{j}) $. 
The distance of $\mathbf f_{m}(x,y)$ from the beginning of $\mathbf \Gamma(\mathbf f)$ is defined 
    referring to a properly-chosen quadratic form  
$ \mathbf f_{0}= \Delta_0 x^2-2 \sqrt{N-\Delta_0} x y+y^2$ hypotetically
 located before $\mathbf f_{1}$. Thus we have
     $d(\mathbf f_{m}, \mathbf f_{0}) = \sum_{j=0}^{m-1}  d(\mathbf f_{j+1}, \mathbf f_{j})$ if $m \leq \tau$. 
     The notion is also extended to index $k \tau \leq m < (k+1) \tau$ by setting
      $d(\mathbf f_{m}, \mathbf f_{0}) =d(\mathbf f_{m \bmod \tau}, \mathbf f_{0}) +kR_{\mathbb F}$.
The distance  $d(\mathbf f_{\tau}, \mathbf f_{0} )$ is exactly equal to $R^*= \ln \mathfrak c_{\tau-1}$, which is the regulator $R_{\mathbb F}$, or three times $R_{\mathbb F}$, 
and the distance  $d(\mathbf f_{\frac{\tau}{2}}, \mathbf f_{0} )$ is exactly equal to 
 $\frac{R^*}{2}$, see  \cite{elia1} for a straightforward proof.  
Now, a celebrated formula of Dirichlet's gives the product 
\begin{equation}
   \label{dirichlet}
  h_{\mathbb F} R_{\mathbb F} =  \frac{\sqrt{D}}{2}  L(1, \chi) = 
      - \sum_{n=1}^{\lfloor \frac{D-1}{2} \rfloor} \left(\frac{D}{n}\right)  \ln\left( \sin \frac{n \pi}{D} \right)
\end{equation} 
where  $ h_{\mathbb F}$ is the class field number, $ L(1, \chi)$ is a Dedekind $L$-function,
 $D=N$ if $N\equiv 1 \bmod 4$ or $D=4N$ otherwise, 
 and character $\chi$ is the Jacobi symbol in this case. If we know $h_{\mathbb F}$ exactly, we know $R^*$ exactly and we can proceed to factorization,
 with complexity $O((\log_2 N)^4)$ \cite{lagarias}, conditioned on the computation of $L(1, \chi)$.
The Dirichlet $L(1,\chi_N)$ function can be efficiently
 evaluated using the following expression for the product $h_{\mathbb F} R_{\mathbb F}$ as a function of $N$ 
\begin{equation}
   \label{eqhR}
  h_{\mathbb F} R_{\mathbb F} = \frac{1}{2} \sum_{x \geq 1}  \left(\frac{D}{x}\right) \left(\frac{\sqrt D}{x} \mbox{erfc}\left(x \sqrt{\frac{\pi}{D}}\right)+ E_1\left(\frac{\pi x^2}{D}\right)\right) ~~ .
\end{equation} 
 where the complementary error function $\mbox{erfc}(x)$, and the exponential integral function $E_1(x)$, can be quickly evaluated.   
Once we know $R^*$, with the Shanks' infrastructural method \cite{shanks}
 or some of its improvements \cite{williams,cohen,schoof0}, we can find $\mathbf f_{\frac{\tau-2}{2}}(x,y)$, thus a factor of $N$.
 The goal is to obtain 
 $\mathbf f_{\frac{\tau-2}{2}}(x,y)$ with as few
  steps as possible.  To this end we can perform
1) giant-steps within $\mathbf \Gamma(\mathbf f)$ which are realized by the Gauss composition law of quadratic forms, followed by a reduction of this form to 
$\mathbf \Gamma(\mathbf f)$, and 
 2) baby-steps moving from one quadratic form to the next in $\mathbf \Gamma(\mathbf f)$.
Two operators $\rho^+$ and $\rho^-$ are further defined \cite[p.259]{cohen} to allow small (baby) steps, precisely
\begin{itemize} 
    \item  $\rho^+$  transforms $\mathbf f_m(x,y)$ into $\mathbf f_{m+1}(x,y)$  in  
$\mathbf \Gamma(\mathbf f)$, and is defined as 
$  \rho^+([a,2b,c]) = [\frac{b_1^2-N}{a},2b_1, a]  ~~, $
where $b_1$ is $2b_1= [2b \bmod (2a)] +2ka$ with $k$ chosen in such a way
 that $-|a| < b_1 < |a|$.
     \item  $\rho^-$  transforms $\mathbf f_m(x,y)$ into $\mathbf f_{m-1}(x,y)$
 in  $\mathbf \Gamma(\mathbf f)$ and  is defined as 
$  \rho^-([a,2b,c]) = [c,2b_1, \frac{b_1^2-N}{c}]  ~~, $
where $b_1$ is $2b_1=[ -2b \bmod (2c)] +2kc$ with $k$ chosen in such a way
 that $-|c| < b_1 < |c|$.
 \end{itemize}

\noindent
 The composed form  $\mathbf f_{m}\bullet \mathbf f_{n}$ has the distance
$d( \mathbf f_{m} \bullet \mathbf f_{n}, \mathbf f_{0}) \approx  d(\mathbf f_{m}, \mathbf f_{0})+ d(\mathbf f_{n}, \mathbf f_{0})$.

\begin{enumerate}
   \item  By the law $\bullet$, $\mathbf \Gamma(\mathbf f)$ resembles a cyclic group, with 
   $\mathbf f_{\tau-1}$ playing the role of identity.

   \item Since in $\mathbf{\Gamma}(f)$ the two triples of signs (signatures) $(-,+,+)$ and $(+,-,-)$ alternate,
  the composed form $\mathbf f_{m}(x,y) \bullet \mathbf f_{n}(x,y)$ must have one of these
    signatures.
\item 
 The composition of a quadratic form with itself is called doubling and denoted 
$2\bullet \mathbf f_{n}$, thus $s$ iterated doublings are written as $2^s\bullet \mathbf f_{n}(x,y)$.
  The distance is nearly maintained by the composition $\bullet$ (giant-steps).
   The  error affecting this distance estimation is of order $O(\ln N)$ as shown by Schoof in \cite{schoof0},
  and is rigorously maintained by the one-step moves $\rho^\pm$ (baby-steps).
\end{enumerate} 

\noindent
An outline of the procedure is the following, assuming that $R^*$ is preliminarily computed:

\begin{enumerate}
    \item Let $\ell$ be a small integer. Compute an initial quadratic form 
       $\mathbf f_\ell =[\Delta_\ell, 2\Omega_\ell, +\Delta_{\ell-1}]$ and 
      its distance $d_\ell=d(\mathbf f_\ell,\mathbf f_0)$ from the continued fraction
       expansion of $\sqrt N$ stopped at term $\ell+1$.
    \item  Compute  $j_t= \lceil \log_2 \frac{R^*}{d_\ell} \rceil$
  \item  Starting with $[\mathbf f_\ell, d_\ell]$, iteratively compute and store in a vector 
  $ \mathcal F_{j_t}$
   the sequence $[ 2^{j} \bullet \mathbf f_\ell, 2^{j}d_\ell]$ up to $j_t$.
 The middle term (i.e. $\mathbf f_{\frac{\tau-2}{2}}$)  of $\mathbf \Gamma(\mathbf f)$  is located between the terms 
       $ 2^{j_t-1} \bullet \mathbf f_\ell$  and  
       $ 2^{j_t}\bullet \mathbf f_\ell$.
  \item The middle term of $\mathbf \Gamma(\mathbf f)$ can be quickly reached using the elements of
     $ \mathcal F_{j_t}$, starting by computing 
$\mathbf f_r =(2^{j_t-1}\bullet  \mathbf f_\ell) \bullet (2^{j_t-2}\bullet  \mathbf f_\ell)$ and checking whether
     $2^{j_t-1}d_\ell+2^{j_t-2} d_\ell$ is greater or smaller than $\frac{R^*}{2}$; in the first case set 
$\mathbf f_s=\mathbf f_r$, otherwise set 
$\mathbf f_s=2^{j_t-1}\bullet  \mathbf f_\ell$. 
  Iterate this composition by computing
          $\mathbf f_r = \mathbf f_s \bullet (2^{i}\bullet  \mathbf f_\ell)$ and setting $\mathbf f_s = \mathbf f_r$
            for decreasing $i$ up to $0$, and let 
          the final term be $[\mathbf f_s, d_s]$.
  \item Iterate the operation $\rho^\pm$  a convenient number $O(\ln N)$ of times, until a
    factor of $4N$  is found.  
\end{enumerate}

\section{Conclusions}
An iterative algorithm has been described which produces a factor of a composite square-free $N$ with $O((\ln (N))^4)$
 iterations at most, if $hR$ is exactly known, $h$ being the class number, and $R$ the regulator of 
 $\mathbb Q(\sqrt N)$. The bound $O((\ln (N))^4)$ is computed by multiplying the number of giant-steps,
 which is $O(\ln (N))$, by the number of steps at each reduction, completing a giant-step, which is upper bounded by $O((\ln (N))^3)$ as shown in \cite{lagarias,schonhage}.
It is remarked that, in this bound computation, the cost of the arithmetics in $\mathbb Z$, i.e.
multiplications and additions of big integers, is not counted  \cite{lagarias}.
Furthermore, it is not difficult to modify the algorithm to use a rough approximation of $hR$;
the computations become cumbersome, but asymptotically the algorithm is polynomial,
because a sufficient approximation of $hR$ is easily obtained by computing the series
 in equation (\ref{eqhR}) truncated at a number of terms $O(\ln (N))$, since the series
converges exponentially  \cite[Proposition 5.6.11, p.262-263]{cohen}.
It remains to ascertain whether this asymptotically-good factoring algorithm is also practically better
than any sub-optimal probabilistic factoring algorithm.


\end{document}